\documentclass[10pt,a4paper]{article}
\usepackage{url}
\usepackage[utf8]{inputenc}
\usepackage[T1]{fontenc}
\usepackage[centertags]{amsmath}
\usepackage{amsfonts}
\usepackage{graphicx}
\usepackage{amsmath}
\usepackage{amssymb}
\usepackage{latexsym}
\usepackage[mathscr]{euscript}
\usepackage{tikz} 
\usetikzlibrary{calc,matrix,arrows,chains,positioning,scopes,decorations.pathmorphing,backgrounds,fit}

\usepackage{amsthm}
\numberwithin{equation}{section}

\theoremstyle{plain}
\newtheorem{teo}{Theorem}[section]
\newtheorem{lem}[teo]{Lemma}
\newtheorem{pro}[teo]{Proposition}
\newtheorem{cor}[teo]{Corollary}

\theoremstyle{definition}
\newtheorem{defi}{Definition}[section]
\newtheorem{notat}[defi]{Notation}
\newtheorem{rem}[defi]{Remark}
\newtheorem{es}[defi]{\textbf{Example}}

\newcommand{\bdfn}{\begin{defi} \begin{rm}}
\newcommand{\edfn}{\end{rm} \end{defi}}
\newcommand{\bthm}{\begin{teo}}
\newcommand{\ethm}{\end{teo}}
\newcommand{\bprop}{\begin{pro}}
\newcommand{\eprop}{\end{pro}}
\newcommand{\bcor}{\begin{cor}}
\newcommand{\ecor}{\end{cor}}
\newcommand{\blem}{\begin{lem}}
\newcommand{\elem}{\end{lem}}
\newcommand{\bfact}{\begin{rem} \begin{rm}}
\newcommand{\efact}{\end{rm} \end{rem}}
\newcommand{\bex}{\begin{es} \begin{rm}}
\newcommand{\eex}{ \end{rm} \end{es}}
\newcommand{\bnot}{\begin{notat} \begin{rm}}
\newcommand{\enot}{\end{rm} \end{notat}}
\newcommand{\ten}{\otimes}
\newcommand{\quot}[2]{{\raisebox{.2em}{$#1$}\left/\raisebox{-.2em}{$#2$}\right.}}

\title{Notes on MV-modules over integral domains}

\author{Serafina Lapenta\\
DiMIE, 
University of Basilicata,\\Viale dell'Ateneo Lucano, 10, \\85100 Potenza,Italy\\serafina.lapenta@unibas.it}
\date{}
\begin{document}

\maketitle
\begin{abstract}
An MV-module is an MV-algebra endowed with a scalar multiplication with scalars in a PMV-algebra (i.e. an MV-algebra endowed with a binary ``ring-like'' product). We investigate the class of semisimple MV-modules over a semisimple and totally ordered integral domain, and prove an adjunction with a special class of linear spaces.
\end{abstract}

\section{Introduction}
MV-algebra were defined in 1958 as the algebraic counterpart \L ukasiewicz infinite-valued logic. They are structures $(A, \oplus, ^*, 0)$ such that $(A, \oplus, 0)$ is an abelian monoid, $x^{**}=x$, and the equations $(x^*\oplus y)^*\oplus y=( y^*\oplus x)^*\oplus x$  and $x \oplus 0^*=0^*$ are satisfied for any $x,y\in A$. The literature on the subject is very wide and we suggest \cite{Cha1, CDM, MunBook} for further details.

The  standard model for an MV-algebra is the unit interval $[0,1]$, with $x\oplus y=\min(x+y,1)$ and $x^*=1-x$ and it generates the variety of MV-algebra. One of the most important achievement in the field is the categorical equivalence between MV-algebras and lattice-ordered groups with strong unit ($\ell u$-groups). We suggest \cite{Birk, Stein} for further details on $\ell u$-groups and related structures. \\
In more details, given an $\ell u$-group $(G,u)$, the interval $[0,u]_G=\{x\in G \mid 0\le x \le u \}$ is an MV-algebra  if we define $x\oplus y=(x+y)\wedge u$ and $x^*=u-x$. This gives us a functor, denoted by $\Gamma$, from $\mathbf{auG}$ (the category whose objects are Abelian lattice-ordered groups with strong unit and whose morphisms are maps that are at the same time groups homomorphisms and lattices homomorphisms) to $\mathbf{MV}$ (the category whose objects are MV-algebras and whose morphisms are homomorphism of MV-algebras).

Product MV-algebras (PMV-algebra for short) are obtained when we endow an MV-algebras with a binary and internal ``ring-like'' product  \cite{MonPMV, DiND}. When the product is a scalar one, with scalars chosen in a PMV-algebras, we obtain the notion of MV-module. A Riesz MV-algebra is an MV-module over $[0,1]$.  We remark that  $[0,1]$ can be seen as a PMV-algebra or a MV-module over itself, when the product (either scalar or internal) coincide with the usual product between real numbers.\\
The functor $\Gamma$ naturally extends to PMV-algebras and MV-modules. We obtain a categorical equivalence between PMV-algebras and a proper subclass of lattice-ordered rings with strong unit ($\ell u$-rings) \cite{DiND}; we will denote this functor by $\Gamma_{(\cdot)}$. In the same way, in \cite{LeuMod} it is proved the categorical equivalence between MV-modules over a fixed PMV-algebra $P$ and lattice-ordered modules over the $\ell u$-rings that corresponds to $P$ via $\Gamma_{(\cdot)}$; we will denote this functor by $\Gamma_{R}$, where $\Gamma_{(\cdot)}(R,u)=P$.

In this short paper we study MV-modules over a special class of PMV-algebras, that is PMV-algebras without zero-divisors. The main result is an adjunction between such MV-modules and linear spaces over a totally ordered and Archimedean field. In order to apply some fundamental results from literature, we need to restrict our work to totally ordered and semisimple PMV-algebras.

\section{Preliminaries}

The main tool in our development is the tensor product of MV-algebras, defined by Mundici in \cite{Mun} and further investigated in \cite{LeuTens, LLTP1, LLTP2}.

Given two MV-algebras $A$ and $B$, the tensor product is the MV-algebra $A\ten_{mv} B$ uniquely defined by the universal bimorphism $\beta\colon A\times B \rightarrow A\ten_{mv}B$ such that $\beta (a,b)=a\ten_{mv}b$. We recall that a bimorphism is a bilinear map that commutes with $\vee$ and $\wedge$ on both argument.

An important subclass of MV-algebras is the one of semisimple algebras. An MV-algebra $A$ is semisimple if the intersection of all maximal ideal (called Radical of $A$, and denoted by Rad($A$)) is zero. We have that, via $\Gamma$, semisimple MV-algebras correspond to Archimedean $\ell u$-groups, where an $\ell$-group $G$ is said to be {\em Archimedean} if $nx\le y$ for any $n\in \mathbb{N}$ and $x\ge 0$, implies $x=0$.

Since the class of semisimple MV-algebras is not closed under tensor product, in \cite{Mun} the semisimple tensor product $\ten_{ss}$ is defined as the quotient
\[ 
A \ten_{ss}B = \quot{A\ten_{mv}B}{Rad(A\ten_{mv}B)},
\]
for any $A$ and $B$ semisimple MV-algebras. It satisfies the following universal property, with respect to semisimple MV-algebras:

for any semisimple MV-algebra $C$ and for any bimorphism $\beta\colon A\times B\rightarrow C$, there is a unique homomorphism of MV-algebras $\omega \colon A\otimes_{ss}B\rightarrow [0, \beta (1,1)]\le_i C$ such that $\omega \circ \beta_{A,B}=\beta$,\\
where $\beta_{A,B}\colon A\times B \rightarrow A\ten_{ss}B$ is defined by $\beta_{A;B}(a,b)=a\ten_{ss}b$.\\
The notation $[0, a]\le_i A$ means that $[0,a]$ is an interval MV-algebra of $A$. See \cite{Mun, LeuTens} for further details.

In \cite{LLTP1} the following is proved.

\bthm \label{teo:TPmod}
Let $A$ be a unital and semisimple PMV-algebra, and $B$ be a semisimple MV-algebra. Then $A\ten_{ss} B$ is an $A$-MV-module.
\ethm
Moreover, denoted by $\mathcal{U}_A(M)$ the MV-reduct of $M$, an MV-module with scalars in $A$, the following universal property holds.

\bthm \cite{LLTP1} \label{un1}
Let $A$ be a  unital, semisimple and totally ordered PMV-algebra, let $B$ be a semisimple MV-algebra. Then for any unital and semisimple $A$-MV-module $M$ and for any homomorphism of MV-algebras $f\colon B\rightarrow \mathcal{U}_{A}(M)$ there is a unique homomorphism of $A$-MV-modules 
$\widetilde{f}:A\ten_{ss} B\rightarrow M$ such that $\widetilde{f}\circ \iota_{B,A}=f$, where $\iota_{B, A} \colon B \rightarrow A\ten_{ss} B$ is the embedding in the tensor product.
\ethm

In \cite{BVanR} the lattice-ordered counterpart of $\ten_{ss}$ is introduced: the authors define the tensor product of Archimedean lattice-ordered groups with strong unit. Given $(G, u_G)$ and $(H, u_H)$ $\ell u$-groups, $(G\ten_a H, u_G \ten_a u_H)$ is an $\ell u$-group  uniquely defined, up to isomorphism, by a universal property with respect to Archimedean structures. 

In \cite{LLTP1} the following is proved.
\bthm \label{isoTPss}
If $(G_A, u_A)$, $(G_B, u_B)$ are Archimedean  $\ell u$-groups and  $A$, $B$ are semisimple MV-algebras such that  $A\simeq \Gamma (G_A, u_A)$ and  $B\simeq \Gamma (G_B, u_B)$ then $A\ten_{ss} B\simeq  \Gamma (G_A \ten_{a}G_B, u_A \ten _{a}u_B)$.
\ethm
\bfact
In \cite{DiND} PMV-algebras are defined in the most general case, while in \cite{MonPMV} any PMV-algebra is unital and commutative. In the sequel we will use the definition from \cite{MonPMV}.
\efact
\section{MV-domains}
We start this section with the definition of an MV-domain.
\bdfn
A PMV-algebra $P$ is called MV-domain if $x\cdot y=0$ implies $x=0$ or $y=0$.
\edfn
\bfact
In \cite{Mon+}, Montagna defines the quasi variety of $PMV^+$-algebras, as PMV-algebras that satisfies the quasi-identity 

$x^2=0$ implies $x=0$.\\
$PMV^+$-algebras are therefore algebras without nilpotent elements, and by definition any MV-domain is a $PMV^+$-algebras. The converse is not true in general.\\
We recall that for $PMV^+$-algebras several important results holds, like the subdirect representation theorem.
\efact

\bprop \label{pro:dom}
Let $P$ be a totally ordered PMV-algebra. $P$ is an MV-domain if and only if the corresponding $\ell u$-ring is a integral domain.
\eprop
\begin{proof}
One direction is obious. For the other direction, let  $P$ be a $MV$-domain such that $P=\Gamma_{(\cdot)}(R,u)$, with $(R,u)$ $\ell u$-ring.\\
Let $x,y$ be elements of $R^+$ such that $x\cdot y=0$. There exist $x_1, \dots x_n, y_1, \dots y_m$ in $P$ such that $x=\sum_{i=1}^n x_i$, $y=\sum_{i=1}^m y_j$. Therefore
\[
x\cdot y= \sum_{i,j}x_i \cdot y_j=0.
\]
Hence for any $i=1,\dots ,n$ and $j=1,\dots ,m$  $x_i \cdot y_j=0$. By hypothesis we have:

(i) there exists one $i$ such that $a_i \neq 0$, then we have all $b_j=0$, and then $b=0$;

(ii) for any $i$ we have $a_i=0$, then $a=0$.\\
The result follows from the total order on $R$.
\end{proof}

\noindent In the follow, we will denote by $\mathbf{MVArDomP}$ the category whose objects are semisimple MV-modules over a semisimple and totally ordered MV-domain $P$, and whose morphisms are homomorphisms of MV-modules.

\bfact \label{rem:MV-dom}
(i) A $P$-ideal $I$ for a $P$-MV-module $A$ is an ideal that satisfies the condition $\alpha x \in I$ for any $\alpha \in P$ and  any $x\in A$. This condition is always satisfied when $P$ is a unital PMV-algebra.

(ii) By \cite[Proposition 3.16]{LeuMod} any object in $\mathbf{MVArDomP}$ is a subdirect product of totally ordered $P$-MV-modules.

(ii) By \cite[Chapter XIV Section 6 Lemma 2]{Birk}, in a totally ordered and Archimedean $\ell$-group any positive element is a strong unit. In particular the product-unit is a strong unit.
\efact

\section{The categorical adjunction}
\noindent Let $\mathbf{LinSpArK}$ be the category whose objects are Archimedean and lattice-ordered linear spaces with strong unit over  $K$,  Archimedean and totally ordered field with strong unit, and whose morphisms are homogeneous homomorphisms of $\ell$-groups.
\bprop \label{pro:01}
Let $(V,u)$ be an object in $\mathbf{LinSpArK}$, and $h\colon V_1\rightarrow V_2$ a morphism between objects $(V_1,u_1)$ and $(V_2, u_2)$ of $\mathbf{LinSpArK}$.  Denoted by $P$ the PMV-algebra $\Gamma_{(\cdot)} (K,e)$, where $e$ is the unit in $K$, $\Gamma_{(K,e)} (V,u)$ is an element of $\mathbf{MVArDomP}$, the category of Archimedean MV-modules over $P$. Moreover, $h|_{\Gamma_{(K,e)} (V_1, u_1)}$ is an homomorphism of MV-modules $\Gamma_{(K,e)} (V_1, u_1)$ and $\Gamma_{(K,e)} (V_2, u_2)$.
\eprop
\begin{proof}
It follows directly from Remark \ref{rem:MV-dom} and \cite[Proposition 4.1]{LeuMod}.
\end{proof}

\blem \label{lem:02}
If $R$ is an Archimedean and totally ordered integral domain, its quotient field $F$ is Archimedean and totally ordered.
\elem
\begin{proof}
$F$ is totally ordered by \cite[Theorem 10.4]{Bly}. Let $a,b \in F^+$ such that $na\le b$ for any $n\in \mathbb N$. By definition, this comes to $n\frac{x_1}{y_1}\le \frac{x_2}{y_2} $, with $x_1,x_1 \in R^+$ and $y_1,y_2 \in R^+ \setminus \{0 \} $ such that $a=\frac{x_1}{y_1},\ b=\frac{x_2}{y_2}$. The latter is equivalent to $\frac{x_2y_1-nx_1y_2}{y_2y_1}\in F^+.$, therefore $x_2y_1-nx_1y_2 \in R^+$ and $nx_1y_2 \leq x_2y_1$.   Since $R$ is an Archimedean integral domain, we get $x_1y_2=0$ and $a=0$. Trivially, the unit in $R$ is unit in $F$.
\end{proof}

\bthm \label{teo:00}
Let $M$ be an object in the category $\mathbf{MVArDomP}$. There exists an Archimedean and lattice-ordered linear space with strong unit $(V,u)$ over a totally ordered and Archimedean field $(K,e)$ uniquely associated to M.
\ethm
\begin{proof} By \cite[Corollary 4.8]{LeuMod}, there exists an Archimedean $\ell u$-group $(G, u)$ and a totally ordered and Archimedean $\ell u$-ring $(R, e)$ such that $P \simeq \Gamma_{(\cdot)}(R,e)$ and $M\simeq \Gamma_{(R,e)}(G,u)$. By \cite[Theorem 3.3]{DiND}, $e$ is unit in $R$ and by Proposition \ref{pro:dom} it is a integral domain. By Lemma \ref{lem:02} the quotient field $K=\{ \frac{a}{b}\ \mid \ a,b\in R\quad b\neq 0\}$ is Archimedean, totally ordered and unital.\\
By Theorem \ref{isoTPss}, $\Gamma (K \ten_aG, e\ten_a u)\simeq \Gamma (K,e)\ten_{ss} \Gamma (G,u)$ and by Theorem \ref{teo:TPmod}, $\Gamma (K,e)\ten_{ss} \Gamma (G,u)$ is a MV-module over $\Gamma (K,e)$, then by \cite[Corollary 4.8]{LeuMod} $K \ten_a G$ is $\ell$-module over $K$ and since $K$ is a field, $K \ten_aG \in \mathbf{LinSpArK}$.\\
The uniqueness of $K \ten_aG $ follows by construction.
\end{proof}

\bprop \label{pro:univ}
Let $M$ be an object in $\mathbf{MVArDomP}$, with $P$ semisimple and totally ordered MV-domain. Let $(G,v)$ be the $\ell u$-group  such that $M\simeq \Gamma (G,v)$, let $(R,e)$ be the integral domain such that $P= \Gamma_{(\cdot)} (R,e)$ and let $K$ be the quotient field of $R$. For any object $(V,u)$ in $\mathbf{LinSpArK}$ and any $f:M\rightarrow \Gamma_{(K,e)} (V,u)$ homomorphism of $P$-MV-modules there exists unique $f^{\sharp}\colon K\ten_a G\rightarrow V$ morphism in $\mathbf{LinSpArK}$ such that $\Gamma_{(K,e)}(f^{\sharp}) \circ \iota_M =f$.
\eprop
\begin{proof}
By definition, $\Gamma_{(K,e)} (V,u)$ is a $\Gamma_{(\cdot)}(K,e)$-MV-module and since $P\subseteq \Gamma(K,e)$, $f$ is well defined as homomorphisms of $P$-MV-modules.\\
By Theorem \ref{un1}, there exists $f^*\colon \Gamma(K,e)\ten_{ss} M \rightarrow \Gamma (V,u)$, homomorphism of $\Gamma(K,e)$-\textit{MV}-modules. By Theorem \ref{isoTPss}, $\Gamma(K,e) \ten_{ss} M \simeq \Gamma (K \ten_a G, e\ten_a v)$. Therefore by \cite[Corollary 4.8]{LeuMod}, $f^*$ extends in a unique way to $f^{\sharp}\colon K \ten_a G \rightarrow V$, morphism in $\mathbf{LinSpArK}$. We remark that by Theorem \ref{un1} $f^*\circ \iota_M=f$, where $\iota_M$ is the standard embedding of $M$ in $\Gamma (K,e)\ten_{ss} M$.
\end{proof} 

\bprop \label{teo:02}
Let $h$ be a morphism between the two objects $M$ and $N$ in the category $\mathbf{MVArDomP}$, with $P\simeq \Gamma_{(\cdot)} (R,e)$, $M\simeq \Gamma (G,v_G)$, $N\simeq \Gamma (H,v_H)$ and let $K$ be the quotient field of $R$. Then there exists a unique morphism $h^{\sharp}\colon K\ten_aG \rightarrow K\overline{\ten  }H$ in $\mathbf{LinSpArK}$ such that $\Gamma_{(K,e)}(h^{\sharp}) \circ \iota_M = \iota_N \circ h$.
\eprop
\begin{proof}
Let $\iota_M$ and $\iota_N$ be the standard embeddings in the tensor products \cite{LLTP1}.By Theorem \ref{isoTPss}, $\Gamma(K,e)\ten_{ss} M\simeq \Gamma (K\ten_aG, e\ten v_G)$ and $\Gamma(K,e)\ten_{ss} N\simeq \Gamma (K\ten_aH, e\ten v_H)$. With abuse of notation, we will denote by $\iota_M$ and $\iota_N$ the composite maps from $M$ and $N$ in $\Gamma (K\ten_aG, e\ten v_G)$ and $\Gamma (K\ten_aH, e\ten v_H)$ rispectively.\\
By Proposition \ref{pro:univ} applied on $\iota_N \circ h$ and $\iota_M$ there exists a unique $h^*\colon \Gamma (K,e)\ten M\rightarrow\Gamma (K\ten_aH, e\ten_av_H)$, such that $h^* \circ \iota_M = \iota_N \circ h$.
\begin{center}
\begin{tikzpicture}
  \node (A) {$ M$};
  \node (B) [right=3cm of A] {$ N$};
  \node (C) [below=2cm of A] {$\Gamma (K,e) \ten M$};
  \node (D) [below=2cm of B] {$\Gamma (K,e) \ten N$};
  \draw[->, above=0.5] (A) to node {$h$} (B);
  \draw[->, left=0.5] (A) to node [swap] {$\iota_{M}$} (C);
  \draw[->, dashed, above=0.5] (C) to node [swap] {$h^*$} (D);
  \draw[->, right=0.5] (B) to node {$\iota_{N}$} (D);
  \node (K) at ($(C)!0.5!(D)$) { };
  \node [below=0.5cm, align=flush center] at (K){ Figure 1};
 \end{tikzpicture}
\end{center}
Again by Theorem \ref{isoTPss} and \cite[Corollary 4.8]{LeuMod} there exists a map $h^{\sharp}\colon K\ten_aG\rightarrow K\ten_aH$. The uniqueness of $h^*$ gives us the desired conclusion.
\end{proof}

\noindent Let $P$ be a semisimple and totally ordered MV-domain, let $(R,e)$ be the $\ell u$-ring such that $P= \Gamma_{(\cdot)} (R,e)$, and let $K$ be the quotient field $R$. We have two functors:
\begin{itemize}
\item $\Gamma _{(K,e)}\colon \mathbf{LinSpArK} \rightarrow \mathbf{MVArModP}$, which is the functor from \cite{LeuMod};
\item $\mathcal{L}\colon \mathbf{MVArModP} \rightarrow \mathbf{LinSpArK}$ such that\\
for any $P$-MV-module $M$, $\mathcal{L}(M)$ is the linear space $K\ten_a G$ defined in Theorem \ref{teo:00},\\
for any morphism $h$, $\mathcal{L}(h)$ is the map $h^{\sharp}$ defined in Proposition \ref{teo:02}. 
\end{itemize} 
\blem
$\mathcal{L}$ is a functor.
\elem
\begin{proof}
Let $h\colon A\rightarrow B$ and $g\colon B\rightarrow C$ be homomorphisms of $P$-MV-modules, with $A=\Gamma (G,u_G)$, $B=\Gamma(H, u_H)$, $C=\Gamma (L, u_L)$.\\
As in Proposition \ref{teo:02}, there exists $h^*\colon \Gamma (K,e)\ten A\rightarrow \Gamma (K\ten_aH, e\ten u_H)$, such that $h^* \circ \iota_A = \iota_B \circ h$, and there exists $g^*\colon \Gamma (K,e)\ten B\rightarrow \Gamma (K\ten_aL, e\ten u_L)$, such that $g^* \circ \iota_B = \iota_C \circ g$. Then 

$ (g^* \circ h^*)\circ \iota_A = g^* \circ (h^*\circ \iota_A)= g^* \circ (\iota_B \circ h)=(g^* \circ \iota_B)\circ h= \iota_C \circ (g\circ h).$
Therefore, $(g\circ h)^*=g^* \circ h^*$. Since $\mathcal{L}(g\circ h)$ is the extension of $(g^* \circ h^*)$ by the inverse functor of $\Gamma_{(K,e)}$, $(g\circ h)^{\sharp}=g^{\sharp}\circ h^{\sharp}$. 
\end{proof}

\bprop
Let $M$ be an element in $\mathbf{MVArDomP}$. Then \\$(\iota_M )_{M\in \mathbf{MVArModP}}$, with $\iota_M \colon M \rightarrow \Gamma(K,e) \ten M$, are a natural transformation between the identity functor on $\mathbf{MVArDomP}$ and the composite functor $\Gamma_{(K,e)} \mathcal{L}$.
\eprop
\begin{proof}
Let $N,L \in \mathbf{MVArModP}$ and let $h\colon N\rightarrow L$ an homomorphism of $P$-MV-modules. We have to prove that $\Gamma_{(K,e)} \mathcal{L}(h)\circ \iota_N = \iota_L \circ h$. This is straightforward, since by definition $\mathcal{L}(h)$ is the extension on linear spaces of $h^*$, then $\Gamma_{(K,e)} \mathcal{L}(h)=h^*$ and the  conclusion follows by Proposition \ref{teo:02}. 
\end{proof}

\bthm
The pair $(\Gamma_{(K,e)}, \mathcal{L})$ is an adjoint pair.
\ethm
\begin{proof}
$\mathcal{L}$ is a left adjoint of $\Gamma_{(K,e)}$ if, for any element $M\in \mathbf{MVArModP}$, any $(V,u)\in \mathbf{LinSpArK}$, and any homomorphism of $P$-MV-module $h\colon M\rightarrow \Gamma (V,u)$ there exists a morphism in $\mathbf{LinSpArK}$ $h^{\sharp}\colon K\ten_aG \rightarrow V$, where $M\simeq \Gamma (G,v)$, such that $\Gamma_{(K,e)} (h^{\sharp})\circ \iota_M=h$. This is proved in Proposition \ref{pro:univ}.
\end{proof}
\bfact
We remark that we cannot have an equivalence between the categories $\mathbf{MVArDomP}$ and $\mathbf{LinSpArK}$. Indeed if $(R,u)=(\mathbb{Z},1)$, $P=\{0,1\}$ and $M=L_3\in \mathbf{MVArModP}$ then $K=\mathbb{Q}$ and $\Gamma (\mathcal{L}(M))=([0,1]\cap \mathbb{Q})\ten L_3 \ncong L_3$.
\efact

\blem \label{lem:03}
Let $P$ be a totally ordered and semisimple MV-domain such that $P=\Gamma_{(\cdot)} (K,e)$, with $K$ totally ordered and Archimedean field. Let $M$ be an semisimple MV-module over $P$. If $\alpha x=0$, then $\alpha=0$ or $x=0$. 
\elem
\begin{proof}
By \cite[Corollary 4.8]{LeuMod}, there exists a semisimple $\ell$-module with strong unit $(V,u)$ over $K$ such that $M=\Gamma_{(K,e)} (V, u)$.  Since $K$ is a field, $(V,u)$ is actually a linear space. The result follows by the remark that the  property holds in any linear space.
\end{proof}

\bprop
Let $P$ be a totally ordered and semisimple MV-domain such that $P=\Gamma_{(\cdot)} (K,e)$, with $K$ totally ordered and Archimedean field. Let $M$ be an Archimedean MV-module over $P$. Then the map 
\[ 
\iota\colon P \rightarrow M, \quad \iota(a)=a1 
\]
is an embedding of MV-algebras.
\eprop
\begin{proof}
By \cite[Lemma 3.11(a)]{LeuMod}, $\iota (0)=0$; by \cite[Definition 3.1]{LeuMod} if $a+b$ is defined, then $\iota (a+b)=(a+b)1=a1+b1=\iota(a) +\iota (b)$ and $\iota$ is linear; by \cite[Lemma 3.11(f)]{LeuMod}, $(a1)^*=a^*1$, then $\iota (a^*)=\iota (a)^*$. Moreover, $a\oplus b = (a\wedge b^*)+b$. Since $P$ is totally ordered, and any linear map is isotone by \cite[Proposition 3.9]{LeuOp}, it follows that $\iota (a\oplus b)=\iota(a)\oplus \iota(b)$. Finally, by Lemma \ref{lem:03}, $\iota$ is injective.
\end{proof}

\end{document}